\documentclass{amsart}

\usepackage{amsmath}
\usepackage{amsfonts}
\usepackage{amssymb,enumerate}
\usepackage{amsthm}
\usepackage[all]{xy}
\usepackage{hyperref}

\newtheorem{lem}{Lemma}[section]
\newtheorem{cor}[lem]{Corollary}
\newtheorem{prop}[lem]{Proposition}
\newtheorem{thm}[lem]{Theorem}

\newtheorem{Defn}[lem]{Definition}
\newtheorem{Ex}[lem]{Example}
\newtheorem{Question}[lem]{Question}
\newtheorem{Property}[lem]{Property}
\newtheorem{Properties}[lem]{Properties}
\newtheorem{Discussion}[lem]{Remark}
\newtheorem{Construction}[lem]{Construction}
\newtheorem{Subprops}{}[lem]
\newtheorem{Para}[lem]{}

\newenvironment{defn}{\begin{Defn}\rm}{\end{Defn}}
\newenvironment{ex}{\begin{Ex}\rm}{\end{Ex}}
\newenvironment{question}{\begin{Question}\rm}{\end{Question}}

\newenvironment{disc}{\begin{Discussion}\rm}{\end{Discussion}}

\newtheorem{intthm}{Theorem}

\newcommand{\comp}[1]{\widehat{#1}}
\newcommand{\ideal}[1]{\mathfrak{#1}}
\newcommand{\m}{\ideal{m}}
\newcommand{\fm}{\ideal{m}}

\newcommand{\cx}{\operatorname{cx}}	
\newcommand{\pdim}{\operatorname{pd}}	
\newcommand{\pd}{\operatorname{pd}}

\newcommand{\ext}{\operatorname{Ext}}

\newcommand{\rank}{\operatorname{rank}}

\newcommand{\id}{\operatorname{id}}

\newcommand{\HH}{\operatorname{H}}
\newcommand{\Hom}{\operatorname{Hom}}	

\newcommand{\coker}{\operatorname{Coker}}

\newcommand{\tor}{\operatorname{Tor}}

\newcommand{\im}{\operatorname{Im}}

\newcommand{\xra}{\xrightarrow}

\newcommand{\shift}{\mathsf{\Sigma}}

\newcommand{\cone}{\operatorname{Cone}}

\newcommand{\Ker}{\operatorname{Ker}}

\newcommand{\wti}{\widetilde}

\newcommand{\cat}[1]{\mathcal{#1}}
\newcommand{\catx}{\cat{X}}
\newcommand{\caty}{\cat{Y}}
\newcommand{\catm}{\cat{M}}
\newcommand{\catv}{\cat{V}}
\newcommand{\catw}{\cat{W}}
\newcommand{\catg}{\cat{G}}
\newcommand{\catp}{\cat{P}}
\newcommand{\catf}{\cat{F}}
\newcommand{\catfc}{\cat{F}'}
\newcommand{\cati}{\cat{I}}
\newcommand{\cata}{\cat{A}}
\newcommand{\catab}{\mathcal{A}b}

\newcommand{\catgi}{\cat{GI}}
\newcommand{\catgp}{\cat{GP}}

\newcommand{\catgic}{\cat{GI}_C}

\newcommand{\catgc}{\cat{G}_C}
\newcommand{\catgpc}{\cat{GP}_C}

\newcommand{\catac}{\cat{A}_C}
\newcommand{\catbc}{\cat{B}_C}

\newcommand{\catic}{\cat{I}_C}
\newcommand{\catpc}{\cat{P}_C}

\newcommand{\opg}{\cat{G}}

\newcommand{\catpd}[1]{\cat{#1}\text{-}\pd}
\newcommand{\xpd}{\catpd{X}}

\newcommand{\gpd}{\catpd{G}}
\newcommand{\catid}[1]{\cat{#1}\text{-}\id}
\newcommand{\yid}{\catid{Y}}

\newcommand{\aext}{\ext_{\cata}}
\newcommand{\ahom}{\Hom_{\cata}}

\newcommand{\finrescat}[1]{\operatorname{res}\comp{\cat{#1}}}
\newcommand{\proprescat}[1]{\operatorname{res}\wti{\cat{#1}}}
\newcommand{\finrescatx}{\finrescat{X}}
\newcommand{\finrescaty}{\finrescat{Y}}

\newcommand{\finrescatw}{\finrescat{W}}

\newcommand{\propcorescatpc}{\operatorname{cores}\wti{\catp_C(R)}}

\newcommand{\proprescatp}{\proprescat{P}}

\newcommand{\proprescatx}{\proprescat{X}}
\newcommand{\proprescaty}{\proprescat{Y}}
\newcommand{\proprescatv}{\proprescat{V}}

\newcommand{\fincorescat}[1]{\operatorname{cores}\comp{\cat{#1}}}
\newcommand{\propcorescat}[1]{\operatorname{cores}\wti{\cat{#1}}}
\newcommand{\fincorescatx}{\fincorescat{X}}

\newcommand{\propcorescati}{\propcorescat{I}}

\newcommand{\fincorescaty}{\fincorescat{Y}}
\newcommand{\fincorescatv}{\fincorescat{V}}
\newcommand{\fincorescatw}{\fincorescat{W}}
\newcommand{\propcorescatx}{\propcorescat{X}}
\newcommand{\propcorescaty}{\propcorescat{Y}}

\newcommand{\propcorescatw}{\propcorescat{W}}

\newcommand{\gw}{\opg(\catw)}

\newcommand{\gnx}[1]{\opg^{#1}(\catx)}
\newcommand{\gx}{\opg(\catx)}

\renewcommand{\geq}{\geqslant}
\renewcommand{\leq}{\leqslant}
\renewcommand{\ker}{\Ker}

\begin{document}

\bibliographystyle{amsplain}

\author{Sean Sather-Wagstaff}

\address{Sean Sather-Wagstaff, Department of Mathematical Sciences, Kent State University,
  Mathematics and Computer Science Building, Summit Street, Kent OH
  44242, USA}
\email{sather@math.kent.edu}
\urladdr{http://www.math.kent.edu/~sather}

\author{Tirdad Sharif}
\address{Tirdad Sharif, School of Mathematics, Institute for Studies in
Theoretical Physics and Mathematics, P. O. Box 19395-5746, Tehran, Iran}
\email{sharif@ipm.ir}
\urladdr{http://www.ipm.ac.ir/IPM/people/personalinfo.jsp?PeopleCode=IP0400060}
\thanks{TS is supported by a grant from IPM,(No. 83130311).}

\author{Diana White} 
\address{Diana White, Department of Mathematics, University of Nebraska,
   203 Avery Hall, Lincoln, NE, 68588-0130 USA} \email{dwhite@math.unl.edu}
\urladdr{http://www.math.unl.edu/~s-dwhite14/}

\title{Stability of Gorenstein Categories }

\date{\today}


\begin{abstract}
We show that an iteration of the procedure used to define the
Gorenstein projective modules over a commutative ring $R$
yields exactly the Gorenstein projective modules.  Specifically,
given
an exact sequence of Gorenstein projective $R$-modules 
$G=\cdots\xra{\partial^G_2}G_1\xra{\partial^G_1}G_0\xra{\partial^G_0} \cdots $
such that  the complexes
$\Hom_R(G,H)$ and $\Hom_R(H,G)$ are exact for each 
Gorenstein projective $R$-module $H$, the module $\coker(\partial^G_1)$
is Gorenstein projective.  The proof of this result hinges upon
our analysis of Gorenstein subcategories of
abelian categories.
\end{abstract}

\maketitle

\section*{Introduction}  

Let $R$ be a commutative ring.
Building
from Auslander and Bridger's work~\cite{auslander:adgeteac,auslander:smt}
on modules of finite G-dimension,
Enochs and Jenda~\cite{enochs:gipm}
and Holm~\cite{holm:ghd} introduce and study the
\emph{Gorenstein projective} $R$-modules as
the modules of the form $\coker(\partial^P_1)$
for some exact sequence of projective $R$-modules 
$$P=\cdots\xra{\partial^P_2}P_1\xra{\partial^P_1}P_0\xra{\partial^P_0} \cdots $$
such that the complex
$\Hom_R(P,Q)$ is exact for each projective $R$-module $Q$.
The class of Gorenstein projective $R$-modules 
is denoted $\catg(\catp(R))$.  

In this paper, we investigate the modules that arise from an iteration
of this construction.  To wit, let $\catg^2(\catp(R))$ denote the 
class of $R$-modules $M$ for which there exists 
an exact sequence of Gorenstein projective $R$-modules 
$$G=\cdots\xra{\partial^G_2}G_1\xra{\partial^G_1}G_0\xra{\partial^G_0} \cdots $$
such that  the complexes
$\Hom_R(G,H)$ and $\Hom_R(H,G)$ are exact for each 
Gorenstein projective $R$-module $H$ and $M\cong\coker(\partial^G_1)$.

One checks readily that there is a containment
$\catg(\catp(R))\subseteq\catg^2(\catp(R))$.  We answer a question from the 
folklore of this subject by verifying that
this containment is always an equality. 
This is a consequence of Corollary~\ref{gxgw02}, as is the dual version for
Gorenstein injective $R$-modules; see Example~\ref{GP}.

\begin{intthm} \label{thma}
If $R$ is a commutative ring, then $\catg(\catp(R))=\catg^2(\catp(R))$.
\end{intthm}

The proof of this result is facilitated by the consideration of a more general 
situation.  Starting with a class of $R$-modules
$\catw$, we consider the associated full subcategory $\catg(\catw)$
of the category of $R$-modules whose objects are
defined as above; see Definition~\ref{defG}.  Section~\ref{sec3}
is devoted to the category-theoretic properties of $\catg(\catw)$,
those needed for the proof
of Theorem~\ref{thma} and others. 
For instance, the next result is contained in Proposition~\ref{exact01} and
Theorem~\ref{admepi01}.

\begin{intthm} \label{thmb}
Assume $\ext^i_R(W,W')=0$ for all $W,W'\in\catw$ and all $i\geq 1$. 
The Gorenstein subcategory
$\catg(\catw)$ is an exact category, and it  is closed under
kernels of epimorphisms (or cokernels of monomorphisms)
if $\catw$ is so.
\end{intthm}

Most of the paper focuses on subcategories of an abelian category $\cata$.
The reader is encouraged to keep certain module-categories
in mind.  Specific examples  are provided in Section~\ref{sec4}, and 
we apply our results to these examples
in Section~\ref{sec5}.

\section{Categories and resolutions}\label{sec1}

Here we set some notation and terminology for use throughout this paper.

\begin{defn} \label{notation01}
In this work
$\cata$ is an abelian category.
We use the term ``subcategory'' for a ``full additive subcategory
that is closed under isomorphisms.''
Write $\catp=\catp(\cata)$ and $\cati=\cati(\cata)$ 
for the  subcategories of projective and injective
objects in $\cata$, respectively.
A subcategory $\catx$ of $\cata$ is 
\emph{exact} if it is  closed under direct summands and extensions.

We fix subcategories $\catx$, $\caty$, $\catw$, and $\catv$  of $\cata$ such that
$\catw\subseteq\catx$
and $\catv \subseteq\caty$.
Write $\catx\perp\caty$ 
if $\aext^{\geq1}(X,Y)=0$ for each object $X$ in $\catx$ and each object $Y$ in $\caty$.
For an object $A$ in $\cata$, write $A\perp\caty$ (resp., $\catx\perp A$)
if $\aext^{\geq1}(A,Y)=0$ for each object  $Y$ in $\caty$
(resp., if $\aext^{\geq1}(X,A)=0$ for each object  $X$ in $\catx$).
We say that $\catw$ is a \emph{cogenerator} for $\catx$ if,
for each object $X$ in $\catx$, there exists an exact sequence in $\catx$
$$0\to X\to W\to X'\to 0$$
such that $W$ is an object in $\catw$.
The subcategory
$\catw$ is an \emph{injective cogenerator} for $\catx$ if 
$\catw$ is a cogenerator for $\catx$ and $\catx\perp\catw$.
We say that 
$\catv$ is a \emph{generator} for $\caty$ if, 
for each object $Y$ in $\caty$, there exists an exact sequence in $\caty$
$$0\to Y'\to V\to Y\to 0$$
such that $V$ is an object in $\catv$.
The subcategory
$\catv$ is a \emph{projective generator} for $\caty$ if 
$\catv$ is a  generator for $\caty$ and
$\catv\perp\caty$.
\end{defn}

\begin{defn} \label{notation07}
An \emph{$\cata$-complex} is a sequence of 
homomorphisms in $\cata$
$$A =\cdots\xra{\partial^A_{n+1}}A_n\xra{\partial^A_n}
A_{n-1}\xra{\partial^A_{n-1}}\cdots$$
such that $ \partial^A_{n-1}\partial^A_{n}=0$ for each integer $n$; the
$n$th \emph{homology object} of $A$ is
$\HH_n(A)=\Ker(\partial^A_{n})/\im(\partial^A_{n+1})$.
We frequently identify objects in $\cata$ with complexes concentrated in degree 0.

Fix an integer $i$.
The $i$th \emph{suspension} of
a complex $A$, denoted $\shift^i A$, is the complex with
$(\shift^i A)_n=A_{n-i}$ and $\partial_n^{\shift^i A}=(-1)^i\partial_{n-i}^A$.
The \emph{hard truncation} $A_{\geq i}$ is the complex
$$A_{\geq i}=\cdots\xra{\partial^A_{i+2}}A_{i+1}\xra{\partial^A_{i+1}}
A_{i}\to 0$$
and the hard truncations $A_{> i}$, $A_{\leq i}$, and $A_{< i}$ are defined similarly.

The complex $A$ is \emph{$\ahom(\catx,-)$-exact} if the complex
$\ahom(X,A)$ is exact for each object $X$ in $\catx$.  
Dually, it is \emph{$\ahom(-,\catx)$-exact} if the complex
$\ahom(A,X)$ is exact for each object $X$ in $\catx$.  
\end{defn}

\begin{defn} \label{notation07a}
Let $A,A'$ be $\cata$-complexes.
The  Hom-complex $\ahom(A,A')$ is the complex of abelian groups defined as
$\ahom(A,A')_n=\prod_p\ahom(A_p,A'_{p+n})$
with $n$th differential $\partial_n^{\ahom(A,A')}$ given by
$\{f_p\}\mapsto \{\partial^{Y'}_{p+n}f_p-(-1)^nf_{n-1}\partial^A_p\}$.
A \emph{morphism} 
is an element of $\ker(\partial_0^{\ahom(A,A')})$
and $\alpha$ is \emph{null-homotopic} if 
it is in $\im(\partial_1^{\ahom(A,A')})$.
Given a second morphism $\alpha'\colon A\to A'$ we say that
$\alpha$ and $\alpha'$ are \emph{homotopic} if the difference
$\alpha-\alpha'$ is null-homotopic.  The morphism $\alpha$ is a
\emph{homotopy equivalence} if there is a morphism
$\beta\colon A'\to A$ such that 
$\beta \alpha $ is homotopic to $\id_{A}$ and
$\alpha\beta$ is homotopic to $\id_{A'}$.
The complex $A$ is \emph{contractible} if the identity morphism
$\id_A$ is null-homotopic.  When $A$ is contractible, it is exact, as is each of
the complexes
$\ahom(A,N)$ and $\ahom(M,A)$ for all objects $M$ and $N$ in $\cata$.

A morphism of complexes $\alpha\colon A\to A'$
induces homomorphisms
$\HH_n(\alpha)\colon\HH_n(A)\to\HH_n(A')$, and $\alpha$ is a
\emph{quasiisomorphism} when each $\HH_n(\alpha)$ is bijective.
The \emph{mapping cone} of $\alpha$ is the complex
$\cone(\alpha)$ defined as
$\cone(\alpha)_n=A'_n\oplus A_{n-1}$
with $n$th differential 
$\partial^{\cone(\alpha)}_n 
= \Bigl(\begin{smallmatrix}\partial_{n}^{A'} & \alpha_{n-1} \\ 0 & -\partial_{n-1}^{A}
\end{smallmatrix} \Bigr)$.
This definition gives a degreewise split exact sequence
$0\to A'\to\cone(\alpha)\to\shift A\to 0$.
Further, the morphism $\alpha$ is a quasiisomorphism if and only if $\cone(\alpha)$ is exact.
Finally, if $\id_A$ is the identity morphism for $A$, then $\cone(\id_A)$ is contractible.
\end{defn}

\begin{defn} \label{notation03}
A complex $X$ is \emph{bounded} if $X_n=0$ for $|n|\gg 0$. When
$X_{-n}=0=\HH_n(X)$ for all $n>0$, the natural morphism
$X\to\HH_0(X)$ is a quasiisomorphism.  In this event, $X$ is an
\emph{$\catx$-resolution} of $M$ if each $X_n$ is an object in $\catx$, and
the following exact sequence is the \emph{augmented
$\catx$-resolution} of $M$ associated to $X$.
$$X^+ = \cdots\xra{\partial^X_{2}}X_1
\xra{\partial^X_{1}}X_0\to M\to 0$$ 
Instead of writing ``$\catp$-resolution'' we will write ``projective resolution.''
The \emph{$\catx$-projective dimension} of $M$ is the quantity
$$\xpd(M)=\inf\{\sup\{n\geq 0\mid X_n\neq 0\}\mid \text{$X$ is an
$\catx$-resolution of $M$}\}.$$ 
The objects of $\catx$-projective dimension 0 are
exactly the objects of $\catx$.
We set
$$\finrescatx=
\text{the subcategory of objects $M$ of $\cata$ with $\xpd(M)<\infty$.}$$
We define \emph{$\caty$-coresolutions} and \emph{$\caty$-injective dimension}
dually.  The \emph{augmented 
$\caty$-coresolution} associated to a $\caty$-coresolution $Y$ is denoted $^+Y$,
and the $\caty$-injective dimension of $M$ is denoted $\yid(M)$.
We set
$$\fincorescaty= 
\text{the subcategory of objects $N$ of $\cata$ with $\yid(N)<\infty$.}$$
\end{defn}

\begin{defn} \label{notation05}
An $\catx$-resolution $X$ is \emph{$\catx$-proper} (or
simply \emph{proper}) if the
the augmented resolution $X^+$ is $\ahom(\catx,-)$-exact.
We set
$$\proprescatx= 
\text{the subcategory of objects of $\cata$ admitting a proper 
$\catx$-resolution.}$$
One checks readily that $\proprescatx$ is additive.
If $X'$ is an object in $\catx$, then the complex $0\to X'\to 0$ is
a proper $\catx$-resolution of $X'$; hence $X'$ is in $\catx$ and so
$\catx \subseteq\proprescatx$.

Projective resolutions are always $\catp$-proper, and so $\cata$ has enough 
projectives if and only if $\proprescatp=\cata$.
If $M$ is an object in 
$\cata$ that admits 
an $\catx$-resolution $X\xra{\simeq}M$ and a projective resolution
$P\xra{\simeq}M$, then there exists a quasiisomorphism 
$P\xra{\simeq} X$.

\emph{Proper coresolutions} are defined dually, and we set
$$\propcorescaty= 
\text{the subcategory of objects of $\cata$ admitting a proper 
$\caty$-coresolution.}$$
Again, $\propcorescaty$ is additive and $\caty \subseteq \propcorescaty$.
Injective coresolutions are always $\cati$-proper, and so $\cata$ has enough 
injectives if and only if $\propcorescati=\cata$. If $N$ is an object in 
$\cata$ that admits 
a $\caty$-coresolution $N\xra{\simeq}Y$ and an injective resolution
$N\xra{\simeq}I$, then there exists a quasiisomorphism 
$Y\xra{\simeq} I$.
\end{defn}

The next lemmata are standard or have standard proofs; 
for~\ref{perp03} see~\cite[pf.~of (2.3)]{auslander:htmcma},
for~\ref{gencat01} see~\cite[pf.~of (2.1)]{auslander:htmcma},
for~\ref{res01} repeatedly apply  Definition~\ref{notation01},
and for the ``Horseshoe Lemma'' \ref{horseshoe01} see~\cite[pf.~of (8.2.1)]{enochs:rha}.

\begin{lem} \label{perp03}
Let $0\to A_1\to A_2\to A_3\to 0$ be an exact sequence in $\cata$.
\begin{enumerate}[\quad\rm(a)]
\item \label{perp03item1}
If $A_3\perp\catw$, then $A_1\perp\catw$ if and only if $A_2\perp\catw$.
If $A_1\perp\catw$  and  $A_2\perp\catw$, 
then $A_3\perp\catw$
if and only if the given sequence is  $\ahom(-,\catw)$ exact.
\item \label{perp03item2}
If $\catv\perp A_1$, then $\catv\perp A_2$ if and only if $\catv\perp A_3$.
If $\catv\perp A_2$  and  $\catv\perp A_3$, 
then $\catv\perp A_1$
if and only if the given sequence is  $\ahom(\catv,-)$ exact. \qed
\end{enumerate}
\end{lem}

\begin{lem} \label{gencat01}
If $\catx\perp\caty$, then  $\catx\perp\finrescaty$ and $\fincorescatx\perp\caty$.
\qed
\end{lem}

\begin{lem} \label{res01}
If $\catw$ is an injective cogenerator for $\catx$, then every object $X$ in $\catx$ admits a
proper $\catw$-coresolution, so  $\catx \subseteq\propcorescatw$.
If $\catv$ is a projective generator for $\caty$, then every object $Y$ in $\caty$ admits a
proper $\catv$-resolution, so  $\caty \subseteq\proprescatv$. \qed
\end{lem}

\begin{lem} \label{horseshoe01}
Let $0\to A'\to A\to A''\to 0$ be an exact sequence in $\cata$.
\begin{enumerate}[\quad\rm(a)]
\item \label{horseshoe01item1}
Assume that $A'$ and $A''$ 
admit proper $\catx$-resolutions $X'\xra{\simeq} A'$ and $X''\xra{\simeq}A''$.
If the given sequence
is $\ahom(\catx,-)$-exact, then $A$ is in $\proprescatx$
with proper $\catx$-resolution $X\xra{\simeq} A$
such that there exists a commutative diagram 
$$
\xymatrix{
0\ar[r] & 
X' \ar[r]^{\left(\begin{smallmatrix}\id_{X'} \\ 0\end{smallmatrix}\right)}\ar[d]_{\simeq} 
& X \ar[r]^{\left(\begin{smallmatrix}0 & \id_{X''} \end{smallmatrix}\right)}\ar[d]_{\simeq} 
& X'' \ar[r]\ar[d]_{\simeq} & 0 \\
0\ar[r] & A' \ar[r] & A \ar[r] & A'' \ar[r] & 0 \\
}
$$
whose top row is  degreewise
split exact 
and
$\partial^X_n
=\Bigl(\begin{smallmatrix}\partial^{X'}_n & f_n \\ 0 & \partial^{X''}_n\end{smallmatrix}\Bigr)$.
\item \label{horseshoe01item2}
Assume that $A'$ and $A''$ 
admit proper $\caty$-coresolutions $A'\xra{\simeq} Y'$ and $A''\xra{\simeq}Y''$.
If the given sequence
is $\ahom(-,\caty)$-exact, then $A$ is in $\propcorescaty$ with
proper $\caty$-coresolution $A\xra{\simeq} Y$
such that there exists a commutative diagram 
$$
\xymatrix{
0\ar[r] & A' \ar[r]\ar[d]_{\simeq} & A \ar[r]\ar[d]_{\simeq} & A'' \ar[r]\ar[d]_{\simeq} & 0 \\
0\ar[r] 
& Y' \ar[r]^{\left(\begin{smallmatrix}\id_{Y'} \\ 0\end{smallmatrix}\right)} 
& Y \ar[r]^{\left(\begin{smallmatrix}0 & \id_{Y''} \end{smallmatrix}\right)} 
& Y'' \ar[r] & 0 \\
}
$$
whose bottom row is  degreewise
split exact 
and
$\partial^Y_n
=\Bigl(\begin{smallmatrix}\partial^{Y'}_n & g_n \\ 0 & \partial^{Y''}_n\end{smallmatrix}\Bigr)$.
\qed
\end{enumerate}
\end{lem}

\section{Technical Results} \label{sec0}

This section consists of three lemmata which
the reader may wish to skip during the first reading.
The first result is  for use in the proof of Lemma~\ref{coresxw01}.

\begin{lem} \label{noproper01}
Let $0\to N'\to N\to N''\to 0$ be an exact sequence in $\cata$.
Assume that $\catw$ is an injective cogenerator for $\catx$ and
$\catv$ is a projective generator for $\caty$.
\begin{enumerate}[\quad\rm(a)]
\item \label{noproper01item1}
Assume that
$N'$ is an object in $\propcorescatx$, 
and $N''$ is an object in $\propcorescatw$.
If $N'\perp\catw$ and $N''\perp\catw$, then $N$ is an object in $\propcorescatw$.
\item \label{noproper01item2}
Assume that
$N''$ is an object in $\proprescaty$, and
$N'$ is an object in $\proprescatv$.
If $\catv\perp N'$ and $\catv\perp N''$, then $N$ is an object in $\proprescatv$.
\end{enumerate}
\end{lem}

\begin{proof}
We prove part~\eqref{noproper01item1}; the proof of
part~\eqref{noproper01item2} is dual.
Let $X'$ be a proper $\catx$-coresolution of $N'$ and 
set $N_i'=\ker(\partial^{X'}_{-i})$ for $i\geq 0$, which yields an exact sequence
\begin{equation} \label{exact05} \tag{$\ast_i$}
0\to N_{i}'\to X'_{-i}\to N_{i+1}'\to 0.
\end{equation}
Each $N'_{i}$ is an object in $\propcorescatx$. 
By induction on $i$,
Lemma~\ref{perp03}\eqref{perp03item1} implies  $N'_{i}\perp\catw$.  

We will construct
exact sequences
\begin{align}
0\to N\to W_0\to N_{1}\to 0 
\label{exact03} 
\tag{$\circledast_0$} \\
0\to N'_{1}\to N_{1}\to N''_{1}\to 0 
\label{exact03'} 
\tag{$\dagger_0$}
\end{align}
such that~\eqref{exact03} is $\ahom(-,\catw)$-exact,
$W_0$ is an object in $\catw$,
and $N''_{1}$ is an object in $\propcorescatw$ such that $N''_{1}\perp\catw$.
Inducting  on $i\geq 0$, this will yield
exact sequences
\begin{align}
0\to N_i\to W_{-i}\to N_{i+1}\to 0 
\label{exact04} 
\tag{$\circledast_i$} \\
0\to N'_{i+1}\to N_{i+1}\to N''_{i+1}\to 0 
\tag{$\dagger_i$}
\end{align}
such that~\eqref{exact04} is $\ahom(-,\catw)$-exact,
$W_i$ is an object in $\catw$,
and $N''_{i+1}$ is an object in $\propcorescatw$ such that $N''_{i+1}\perp\catw$.
Splicing together the sequences~\eqref{exact04} 
will then yield a proper $\catw$-coresolution of $N_0=N$.

Consider the following pushout diagram
whose top row is the given exact sequence and whose leftmost column is ($\ast_{0}$). 
\begin{equation}
\begin{split} \label{diag05}
\xymatrix{
& 0 \ar[d] & 0 \ar[d] \\
0\ar[r]  & N'\ar[r] \ar[d] \ar@{}[rd]|>>{\lrcorner} & N\ar[r] \ar[d] & N''\ar[r]\ar[d]_{\cong}  & 0 \\
0\ar[r] & X'_{0}\ar[r]\ar[d] & V\ar[r]\ar[d] & N''\ar[r] & 0 \\
& N'_{1}\ar[r]^{\cong}\ar[d] & N'_{1}\ar[d] \\
& 0 & 0
}
\end{split}
\end{equation}
Since $N''\perp\catw$ and $X'_{0}\perp\catw$, the middle row of~\eqref{diag05}
is $\ahom(-,\catw)$-exact and $V\perp\catw$ by Lemma~\ref{perp03}\eqref{perp03item1}.
Note that $X_{0}'$ is in $\propcorescatw$ by Lemma~\ref{res01}.  Since $N''$
is in $\propcorescatw$,
Lemma~\ref{horseshoe01}\eqref{horseshoe01item2} implies that
$V$ is in $\propcorescatw$.
Hence, 
a proper $\catw$-coresolution of $V$ provides an exact sequence 
$$0\to V\to W_{0}\to N''_{1}\to 0$$
that is $\ahom(-,\catw)$-exact
with  objects $W_{0}$ in $\catw$ and $N''_{1}$ in $\propcorescatw$.  
By assumption, we have $\catx\perp\catw$ and $\catw\subseteq\catx$,
and so $W_0\perp\catw$ because $W_0$ is in $\catw$.
Since $V\perp\catw$, Lemma~\ref{perp03}\eqref{perp03item1} 
implies $N_{1}''\perp\catw$.
With the center column of~\eqref{diag05} this yields another pushout diagram, and we shall show that 
the middle row and the rightmost column satisfy the conditions for~\eqref{exact03} 
and~\eqref{exact03'}, respectively.
\begin{equation}
\begin{split} \label{diag06}
\xymatrix{
& & 0 \ar[d] & 0 \ar[d] \\
0\ar[r]  & N\ar[r] \ar[d]_{\cong} & V\ar[r] \ar[d] \ar@{}[rd]|>>{\lrcorner} & N'_{1}\ar[r]\ar[d]  & 0 \\
0\ar[r] & N\ar[r] & W_{0}\ar[r]\ar[d] & N_{1}\ar[r]\ar[d] & 0 \\
& & N''_{1}\ar[r]^{\cong}\ar[d] & N''_{1}\ar[d] \\
& & 0 & 0
}
\end{split}
\end{equation}
We have already seen that
$W_0$ is  in $\catw$ 
and $N''_{1}$ is  in $\propcorescatw$.
Since $N''_{1}\perp\catw$
and $N'_{1}\perp\catw$,
the rightmost column of~\eqref{diag06} with
Lemma~\ref{perp03}\eqref{perp03item1}
yields 
$N_{1} \perp\catw$, and so the middle row of~\eqref{diag06} is
$\ahom(-,\catw)$-exact.
\end{proof}

Next is a key lemma  for  both Theorems~\ref{thma} and~\ref{thmb}
from the introduction.

\begin{lem} \label{coresxw01}
Assume that $\catw$ is an injective cogenerator for $\catx$ and
$\catv$ is a projective generator for $\caty$.
\begin{enumerate}[\quad\rm(a)]
\item \label{coresxw01item1}
If $M$ is an object in $\propcorescatx$ and $M\perp\catw$, then
$M$ is in $\propcorescatw$.
\item \label{coresxw01item2}
If $N$ is an object in $\proprescaty$ and $\catv\perp N$, then
$N$ is in $\proprescatv$.
\end{enumerate}
\end{lem}

\begin{proof}
We prove part~\eqref{coresxw01item1}; the proof of
part~\eqref{coresxw01item2} is dual.
Let $M\xra{\simeq} X$ be a proper $\catx$-coresolution of $M$.  Setting
$M'=\im(\partial^X_0)$ yields an exact sequence 
\begin{equation}\label{exact06}
0\to M\to X_0\to M'\to 0
\end{equation}
that is $\ahom(-,\catx)$-exact.  Since 
$\catw$ is a cogenerator for $\catx$
there is an exact sequence with objects $W_0$ in $\catw$ and $X'$ in $\catx$.  
\begin{equation}\label{exact07}
0\to X_0\to W_0\to X'\to 0
\end{equation}
Consider the pushout diagram
whose top row is~\eqref{exact06} and whose middle column is~\eqref{exact07}.
\begin{equation}
\begin{split} \label{diag07}
\xymatrix{
& & 0 \ar[d] & 0 \ar[d] \\
0\ar[r]  & M\ar[r] \ar[d]_{\cong} & X_0\ar[r] \ar[d] \ar@{}[rd]|>>{\lrcorner} & M'\ar[r]\ar[d]  & 0 \\
0\ar[r] & M\ar[r] & W_{0}\ar[r]\ar[d] & U\ar[r]\ar[d] & 0 \\
& & X'\ar[r]^{\cong}\ar[d] & X'\ar[d] \\
& & 0 & 0
}
\end{split}
\end{equation}
We will show that $U$ is an object in $\propcorescatw$ and that the middle row
of~\eqref{diag07} is $\ahom(-,\catw)$-exact.  It will then follow that a proper
$\catw$-coresolution for $M$ can be obtained by splicing the middle row
of~\eqref{diag07} with a proper $\catw$-coresolution of $U$.

The object $M'$ is in $\propcorescatx$ by construction, and $X'$ and $X_0$ are in $\catx$.  
Thus, $X'$ is an object in $\propcorescatw$ by
Lemma~\ref{res01}, and $X'\perp\catw$ and $X_0\perp\catw$
by hypothesis.
Since 
the top row of~\eqref{diag07} is $\ahom(-,\catx)$-exact,
the assumption
$\catw\subseteq\catx$ implies that the top row of~\eqref{diag07} is $\ahom(-,\catw)$-exact.
Hence,
the assumption $M\perp\catw$ yields $M'\perp\catw$ by Lemma~\ref{perp03}\eqref{perp03item1}.
With the rightmost column of~\eqref{diag07}, Lemma~\ref{noproper01}\eqref{noproper01item1}
implies that $U$ is an object in $\propcorescatw$.
Since $X'\perp\catw$ and $M'\perp\catw$, Lemma~\ref{perp03}\eqref{perp03item1}
yields $U\perp\catw$,
and so  the middle row
of~\eqref{diag07} is $\ahom(-,\catw)$-exact. 
\end{proof}

The last result in this section is a tool for Proposition~\ref{wcogengx01}.

\begin{lem} \label{cogen02}
For $n=0,1,2,\ldots,t$, let $\catx_n$ and $\caty_n$ be subcategories of $\cata$.
\begin{enumerate}[\quad\rm(a)]
\item \label{cogen02item1}
Assume that 
$\catx_n$ is a cogenerator for $\catx_{n+1}$ for each $n\geq 0$ and $\catx_t\perp\catx_0$.
If $\catx_t$ is closed under extensions,
then $\catx_0$ is an injective cogenerator for $\catx_t$. 
\item \label{cogen02item2}
Assume that 
$\caty_n$ is a generator for $\caty_{n+1}$ for each $n\geq 0$ and $\caty_0\perp\caty_t$.
If $\caty_t$ is closed under extensions, then $\caty_0$ is a projective generator for $\caty_t$.
\end{enumerate}
\end{lem}

\begin{proof}
We prove part~\eqref{cogen02item1};  the proof of part~\eqref{cogen02item2} is dual.
Since $\catx_t\perp\catx_0$ by assumption,
it remains to show that $\catx_0$ is a cogenerator for $\catx_t$.
Fix an object $X_t$ in $\catx_t$.  By reverse induction on $i<t$, we will construct
exact sequences
\begin{equation} \label{exact11} \tag{$\ast_i$}
0\to X_t\to X_i\to X_t^{(i)}\to 0
\end{equation}
with objects $X_i$ in $\catx_i$ and $X_t^{(i)}$ in $\catx_t$.
Since $\catx_{t-1}$ is a cogenerator for $\catx_t$, the sequence ($\ast_{t-1}$) 
is known to exist.  By induction, we assume that ($\ast_{i}$) has been constructed
and construct ($\ast_{i-1}$) from it.
From ($\ast_{i}$) we have the object $X_{i}$ in $\catx_{i}$.
Since $\catx_{i-1}$ is a cogenerator for $\catx_i$, there is an exact sequence
\begin{equation} \label{exact10} \tag{$\circledast_i$}
0\to X_i\to X_{i-1}\to X_i'\to 0
\end{equation}
with $X_{i-1}$ in $\catx_{i-1}$ and $X_{i}'$ in $\catx_{i}$.
Consider the pushout diagram
\begin{equation*} 
\xymatrix{
& 0 \ar[d] & 0 \ar[d] \\
& X_t\ar[d] \ar[r]^{\cong} & X_t\ar[d] \\
0 \ar[r] & X_i \ar[r]\ar[d] \ar@{}[rd]|>>{\lrcorner} & X_{i-1}\ar[r]\ar[d] & X_i' \ar[r]\ar[d]_{\cong} & 0 \\
0 \ar[r] & X_{t}^{(i)} \ar[r]\ar[d] & X_{t}^{(i-1)}\ar[r]\ar[d] & X_i' \ar[r] & 0 \\
& 0 & 0
}
\end{equation*}
whose leftmost column is ($\ast_i$) and whose middle row is ($\circledast_i$).
The object $X_{i}'$ is in $\catx_{i}$, and hence in $\catx_{t}$.  Since
$X_{t}^{(i)}$ is also in $\catx_{t}$, the exactness of the bottom row,
with the fact that $\catx_{t}$ is closed under extensions, implies that
$X_{t}^{(i-1)}$ is in $\catx_t$, so the center column of the diagram is the desired sequence
($\ast_{i-1}$).
\end{proof}

\section{Categories of Interest} \label{sec4}

Much of the motivation for this work comes from module categories.
In reading this paper, the reader may find it helpful to keep in mind the examples outlined in the 
next few paragraphs, wherein $R$ is a commutative ring.  

\begin{defn} \label{notation08}
Let $\catm(R)$ denote the category of $R$-modules.
To be clear, we write $\catp(R)$  for the subcategory of projective $R$-modules
and $\cati(R)$ for the subcategory of injective $R$-modules.
If $\catx(R)$ is a subcategory of $\catm(R)$, then $\catx^f(R)$ is the  subcategory
of finitely generated modules in $\catx(R)$.  
Also set $\catab=\catm(\mathbb{Z})$, the category
of abelian groups.
\end{defn}

The study of semidualizing modules was initiated independently (with different names)
by Foxby~\cite{foxby:gmarm}, Golod~\cite{golod:gdagpi},
and Vasconcelos~\cite{vasconcelos:dtmc}.

\begin{defn} \label{notation08a}
An $R$-module $C$ is \emph{semidualizing} if it satisfies the following.
\begin{enumerate}[\quad(1)]
\item $C$ admits a (possibly unbounded) resolution by finite rank free $R$-modules.
\item The natural homothety map $R\to\Hom_R(C,C)$ is an isomorphism.
\item $\ext_R^{\geq 1}(C,C)=0$.
\end{enumerate}
A finitely generated projective  $R$-module of rank 1 is semidualizing.
If $R$ is Cohen-Macaulay, then $C$  is \emph{dualizing}
if it is semidualizing and $\id_R(C)$ is finite.
\end{defn}

Based on the work of Enochs and Jenda~\cite{enochs:gipm},
the following notions were introduced and studied in this generality by
Holm and J\o rgensen~\cite{holm:smarghd}
and White~\cite{white:gpdrsm}.

\begin{defn} \label{notation08b}
Let $C$ be a semidualizing $R$-module, and set
\begin{align*}
\catpc(R)&=\text{the subcategory of modules $P\otimes_R C$ where $P$ is $R$-projective}\\
\catic(R)&=\text{the subcategory of modules $\Hom_R(C,I)$ where $I$ is $R$-injective.}
\end{align*}
Modules in $\catpc(R)$ and $\catic(R)$ are called \emph{$C$-projective}
and \emph{$C$-injective}, respectively.
A \emph{complete $\catp\catpc$-resolution} is a complex $X$ of $R$-modules 
satisfying the following.
\begin{enumerate}[\quad(1)]
\item $X$ is exact and $\Hom_R(-,\catpc(R))$-exact.
\item $X_i$ is projective if  $i\geq 0$ and $X_i$ is  $C$-projective if $i< 0$.
\end{enumerate}
An $R$-module $M$ is \emph{$\text{G}_C$-projective} if there
exists a complete $\catp\catpc$-resolution $X$ such that $M\cong\coker(\partial^X_1)$,
in which case $X$ is a \emph{complete $\catp\catpc$-resolution of $M$}.  We set
$$\catgpc(R)=\text{the subcategory of $\text{G}_C$-projective $R$-modules}.$$
Projective $R$-modules and $C$-projective $R$-modules are $\text{G}_C$-projective,
and $\catpc(R)$ is an injective cogenerator for 
$\catgpc(R)$ by~\cite[(2.5),(2.13)]{holm:smarghd} 
and~\cite[(3.2),(3.9)]{white:gpdrsm}. 

A \emph{complete $\catic\cati$-coresolution}
is a complex $Y$ of $R$-modules 
such that:
\begin{enumerate}[\quad(1)]
\item $Y$ is exact and $\Hom_R(\catic(R),-)$-exact.
\item $Y_i$ is injective if  $i\leq 0$ and $Y_i$ is  $C$-injective if $i>0$.
\end{enumerate}
An $R$-module $N$ is \emph{$\text{G}_C$-injective} if there
exists a complete $\catic\cati$-coresolution $Y$ such that $N\cong\Ker(\partial^Y_0)$,
in which case $Y$ is a \emph{complete $\catic\cati$-coresolution of $N$}.
We set
$$\catgic(R)=\text{the subcategory of $\text{G}_C$-injective $R$-modules}.$$
An $R$-module that is injective or $C$-injective is $\text{G}_C$-injective,
and $\catic(R)$ is a projective generator for 
$\catgic(R)$ by~\cite[(2.6),(2.13)]{holm:smarghd} 
and results dual to~\cite[(3.2),(3.9)]{white:gpdrsm}.
\end{defn}

The next definition was first introduced by Auslander and 
Bridger~\cite{auslander:adgeteac,auslander:smt}
in the case $C=R$, and in this generality
by Golod~\cite{golod:gdagpi} and Vasconcelos~\cite{vasconcelos:dtmc}.

\begin{defn} \label{notation08c}
Assume that $R$ is noetherian  and $C$ is a
semidualizing $R$-module.  A finitely generated $R$-module
$H$ is \emph{totally $C$-reflexive} if 
\begin{enumerate}[\quad(1)]
\item $\ext_R^{\geq 1}(H,C)=0=\ext_R^{\geq 1}(\Hom_R(H,C),C)$, and 
\item the natural biduality map $H\to\Hom_R(\Hom_R(H,C),C)$ is an isomorphism.
\end{enumerate}
Each finitely generated
$R$-module that is either projective or $C$-projective is totally $C$-reflexive.
We set
$$\catgc(R)=\text{the subcategory of totally $C$-reflexive $R$-modules}$$
and $\catg(R)=\catg_R(R)$.
The equality $\catgc(R)=\catgpc^f(R)$ is
by~\cite[(5.4)]{white:gpdrsm}, and $\catpc^f(R)$ is an injective cogenerator for 
$\catgc(R)$ by~\cite[(3.9),(5.3),(5.4)]{white:gpdrsm}.  
\end{defn}

Over a noetherian ring, the next categories were introduced by 
Avramov and Foxby~\cite{avramov:rhafgd}
when $C$ is dualizing, and 
by Christensen~\cite{christensen:scatac} for arbitrary $C$.
(Note that these works (and others) use the notation $\catac(R)$
and $\catbc(R)$ for certain categories of complexes, while our 
categories consist precisely of the modules in these other categories.)
In the non-noetherian setting, these definitions are from~\cite{holm:fear,white:gpdrsm}.

\begin{defn} \label{notation08d}
Let $C$ be a
semidualizing $R$-module.  
The \emph{Auslander class} of $C$ is the subcategory $\catac(R)$
of $R$-modules $M$ such that 
\begin{enumerate}[\quad(1)]
\item $\tor^R_{\geq 1}(C,M)=0=\ext_R^{\geq 1}(C,C\otimes_R M)$, and
\item The natural map $M\to\Hom_R(C,C\otimes_R M)$ is an isomorphism.
\end{enumerate}
The \emph{Bass class} of $C$ is the subcategory $\catbc(R)$
of $R$-modules $N$ such that 
\begin{enumerate}[\quad(1)]
\item $\ext_R^{\geq 1}(C,M)=0=\tor^R_{\geq 1}(C,\Hom_R(C,M))$, and 
\item The natural evaluation map $C\otimes_R\Hom_R(C,N)\to N$ is an isomorphism.
\end{enumerate}
\end{defn}

For a discussion of the next subcategory, consult~\cite[Sec.~5.3]{enochs:rha}.

\begin{defn} \label{fc}
The category of \emph{flat cotorsion} $R$-modules is the subcategory $\catfc(R)$
of flat $R$-modules $F$ such that $\ext^{\geq 1}_R(F',F)=0$ for each flat $R$-module $F'$.
\end{defn}

Gerko~\cite{gerko:ohd} introduced
our final subcategory of interest.

\begin{defn} \label{pci}
Assume that $(R,\m,k)$ is local and noetherian.  The 
\emph{complexity} of a finitely generated $R$-module $M$ is
$$\cx_R(M)
=\inf\{d\in\mathbb{N}\mid \text{there exists $c>0$ such that 
$\beta^R_n(M)\leq cn^{d-1}$ for  $n\gg 0$}\}$$
where $\beta_n^R(M)=\rank_k(\tor_i^R(M,k))$ is the $n$th \emph{Betti number} of $M$.
Let $\catg'(R)$ denote the subcategory
of modules in $\catg(R)$ with finite complexity.
\end{defn}

\section{Gorenstein Subcategories}\label{sec3}

In this section, we introduce and study the Gorenstein subcategory $\catg(\catx)$.

\begin{defn} \label{defG}
An exact complex   in $\catx$ is \emph{totally $\catx$-acyclic} if it is 
$\ahom(\catx,-)$-exact and $\ahom(-,\catx)$-exact.
Let
$\opg(\catx)$ denote the  subcategory of $\cata$ whose objects are
of the form $M\cong\coker(\partial_1^X)$ for some totally $\catx$-acyclic complex $X$;
we say that $X$ is a \emph{complete $\catx$-resolution of $M$}.
Note that the isomorphisms
\begin{align*}
\ahom(X'\oplus X'',X)\cong\ahom(X',X)\oplus\ahom(X'',X) \\
\ahom(X,X'\oplus X'')\cong\ahom(X,X')\oplus\ahom(X,X'')
\end{align*}
show that the direct sum of totally $\catx$-acyclic $\catx$-complexes
is totally $\catx$-acyclic, and hence $\catg(\catx)$ is additive.
Set $\opg^0(\catx)=\catx$ and $\opg^1(\catx)=\opg(\catx)$, and inductively
set $\opg^{n+1}(\catx)=\opg(\opg^n(\catx))$ for $n\geq 1$.
\end{defn}

\begin{disc} \label{gtriv01}
Any contractible $\catx$-complex is totally $\catx$-acyclic; see~\ref{notation07a}.
In particular, for any object $X$ in $\catx$, the complex
$$0\to X\xra{\id_X}X\to 0$$
is a complete $\catx$-resolution, and so 
$X$ is an object in $\opg(\catx)$.
Hence, $\catx \subseteq\opg(\catx)$, and inductively
$\opg^n(\catx) \subseteq\opg^{n+1}(\catx)$
for each $n\geq 0$.

There is a containment  $\catg(\catx)\subseteq\proprescatx\cap\propcorescatx$.
Indeed, If $M$ is an object in $\opg(\catx)$ with complete $\catx$-resolution $X$,
then the hard truncation $X_{\geq 0}$ is a proper $\catx$-resolution of $M$ and
$X_{< 0}$ is a proper $\catx$-coresolution of $M$.
\end{disc}

The orthogonality properties documented next will be very useful in the sequel;
compare to~\cite[(4.2.5)]{christensen:gd}.

\begin{prop} \label{gperp01}
If $\catx\perp\catw$ and $\catv\perp\caty$,
then $\opg^n(\catx)\perp\finrescatw$ and $\fincorescatv \perp \opg^n(\caty)$
for each $n\geq 1$.
In particular, 
if $\catw\perp\catw$, then
$\opg^n(\catw)\perp\finrescatw$ and $\fincorescatw \perp \opg^n(\catw)$
for each $n\geq 1$.
\end{prop}

\noindent \emph{Proof.}
Assuming $\catx\perp\catw$, we will show  $\opg(\catx)\perp\catw$;
the conclusion $\opg^n(\catx)\perp\catw$ will follow by induction, and
$\opg^n(\catx)\perp\finrescatw$ will then follow from Lemma~\ref{gencat01}.
The other conclusion is verified dually.
Let $M$ be an object in $\opg(\catx)$ with
complete $\catx$-resolution $X$, and let $W$ be an object in $\catw$.  
For each integer $i$ set $M_i=\coker(\partial_i^X)$. 
Note that $M\cong M_1$.  
The exact sequence
$$0\to M_{i+1}\xra{\epsilon_i} X_{i-1}\to M_i\to 0$$
is $\ahom(-,\catx)$-exact, and so it is $\ahom(-,W)$-exact.  
In particular, the map $\ahom(\epsilon_i,W)$ is surjective.
Since $\catx\perp\catw$, part of the beginning of the 
associated long exact sequence in $\aext(-,W)$ is
$$\ahom(X_{i-1},W)\xra{\ahom(\epsilon_i,W)}\ahom(M_{i+1},W)
\to \aext^1(M_i,W)\to 0$$
so the surjectivity of $\ahom(\epsilon_i,W)$ implies $\aext^1(M_i,W)= 0$.
In particular, 
$$\aext^1(M,W)\cong\aext^1(M_1,W)=0.$$
If $j\geq 2$, then the remainder of the long exact sequence yields
isomorphisms $\aext^j(M_i,W)\cong\aext^{j-1}(M_{i+1},W)$.
Inductively, this yields the second isomorphism
in the next sequence and the desired vanishing
\begin{xxalignat}{3}
  &{\hphantom{\square}}& \aext^j(M,W)\cong\aext^j(M_1,W)
  &\cong\aext^1(M_{j},W)=0.  &&\square
\end{xxalignat}

We next present a ``Horseshoe Lemma'' for complete $\catx$-resolutions;
compare to~\cite[(4.3.5.a)]{christensen:gd}.

\begin{prop} \label{exten01}
Consider an exact sequence in $\cata$
$$0\to M'\to M\to M''\to 0$$
that is $\ahom(\catx,-)$-exact and $\ahom(-,\catx)$-exact.
If $M'$ and $M''$ are objects in $\opg(X)$, then so is $M$.
Furthermore, given complete $\catx$-resolutions
$X'$ and $X''$ of $M'$ and $M''$, respectively, there is
a degreewise split exact sequence of complexes
$$0\to X'\to X\to X''\to 0$$
such that $X$ is a complete $\catx$-resolution of $M$,
the induced sequence 
$$0\to\coker(\partial^{X'}_1)\to\coker(\partial^{X}_1)\to\coker(\partial^{X''}_1)\to 0$$
is equivalent to the original sequence, and 
$\partial^X_n=\Bigl(\begin{smallmatrix} \partial^{X'}_n & f_n \\ 0 & \partial^{X''}_n
\end{smallmatrix}\Bigr)$ for each $n\in\mathbb{Z}$.
\end{prop}

\begin{proof}
Let $X'$ and $X''$ be complete $\catx$-resolutions for $M'$ and $M''$, respectively.
Lemma~\ref{horseshoe01} yields a degreewise split exact sequence
of complexes $0\to X'\to X\to X''\to 0$ such that
$M\cong\coker(\partial^X_1)$.  Note that each $X_i\cong X_i'\oplus X_i''$ is 
in $\catx$.  Since 
the complexes $X'$ and $X''$ are both $\ahom(\catx,-)$-exact
and $\ahom(\catx,-)$-exact, the same is true of  $X$.  
So, $X$ is a complete $\catx$-resolution of $M$.
\end{proof}

\begin{cor} \label{exten02}
If  $\catw\perp\catw$,
then $\opg(\catw)$ is  closed under extensions.
\end{cor}

\begin{proof}
Proposition~\ref{gperp01} implies 
$\opg(\catw)\perp\catw$ and $\catw \perp \opg(\catw)$.  
Hence, any exact sequence 
$0\to M'\to M\to M''\to 0$ with $M'$ and $M''$ objects 
in $\opg(\catw)$ is
$\ahom(\catw,-)$-exact and $\ahom(-,\catw)$-exact.
Now apply Proposition~\ref{exten01}.
\end{proof}

It is unclear in general
whether $\catw$ is an injective cogenerator for $\opg^n(\catx)$ without 
the extra hypotheses in our next result.

\begin{prop} \label{wcogengx01}
Fix an integer $n\geq 1$.
\begin{enumerate}[\quad\rm(a)]
\item \label{wcogengx01item1}
If  $\catw$ is an injective cogenerator for $\catx$ and 
$\opg^n(\catx)$ is closed under extensions, then $\catw$ is
an injective cogenerator for $\opg^n(\catx)$.
\item \label{wcogengx01item2}
If $\catv$ is a projective generator for $\caty$
and $\opg^n(\caty)$ is closed under extensions, then $\catv$ is
a projective generator for $\opg^n(\caty)$.
\end{enumerate}
\end{prop}

\begin{proof}
We prove part~\eqref{wcogengx01item1}; the proof of part~\eqref{wcogengx01item2}
is dual.  Set $\catx_0=\catw$, $\catx_1=\catx$, and $\catx_t=\opg^{t-1}(\catx)$ for
$t=2,\ldots,n+1$.
By Proposition~\ref{gperp01} we know $\opg^{n}(\catx)\perp\catw$,
so the desired conclusion follows from Lemma~\ref{cogen02}.
\end{proof}

In Section~\ref{sec5} we document the consequences of the following result
for the examples of Section~\ref{sec4}.

\begin{cor} \label{wcogengw01}
If $\catw\perp\catw$, then $\catw$ is both an injective cogenerator and
a projective generator for $\opg(\catw)$.
\end{cor}

\begin{proof}
This follows from Corollary~\ref{exten02} and Proposition~\ref{wcogengx01}.
\end{proof}

The next result extends part of Remark~\ref{gtriv01} and represents
a first step in the proof of Theorem~\ref{thma} from the introduction.

\begin{thm} \label{ggisg01}
Assume that $\catw$ is an injective cogenerator for $\catx$ and that
$\catv$ is a projective generator for $\caty$.  
\begin{enumerate}[\quad\rm(a)]
\item \label{ggisg01item1}
If $\catx$ is closed under extensions, then 
$\opg^n(\catx) \subseteq\propcorescatw$ for each $n\geq 0$.
\item \label{ggisg01item2}
If $\caty$ is closed under extensions, then 
$\opg^n(\caty) \subseteq\proprescatv$ for each $n\geq 0$.
\end{enumerate}
\end{thm}

\begin{proof}
We prove part~\eqref{ggisg01item1} by induction on $n$; the proof of
part~\eqref{ggisg01item2} is dual.
The case $n=0$ is in Lemma~\ref{res01}.  
When $n=1$, note that an object $M$ in $\opg(\catx)$ is
in $\propcorescatx$ by Remark~\ref{gtriv01}, and one has $M\perp\catw$ by 
Proposition~\ref{gperp01}; now apply Lemma~\ref{coresxw01}.

Assume $n>1$ and $\opg^{n-1}(\catx) \subseteq\propcorescatw$.
Fix an object $M$ in $\opg^{n}(\catx)$ and a complete
$\opg^{n-1}(\catx)$-resolution $G$ of $M$.
By definition, the complex $G$ is exact and 
there is an isomorphism $M\cong\ker(\partial^G_{-1})$.
For each integer $j$, set $M_j=\ker(\partial^G_{j})$
and observe that each $M_j$ is an object in $\opg^{n}(\catx)$.
Since each object $G_j$ is in $\opg^{n-1}(\catx)$, Proposition~\ref{gperp01}
implies $M_j\perp\catw$ and $G_j\perp\catw$ for each integer $j$,
and we
consider the exact sequences 
\begin{equation} \label{exact08} \tag{$\maltese_j$}
0\to M_j\to G_j\to M_{j-1}\to 0.
\end{equation}
Our  induction assumption implies that each object $G_j$ is in $\propcorescatw$.

By induction on $i\geq 0$, we construct exact sequences in $\cata$
\begin{gather}
\tag{$\dag_i$}
0\to M\to W_0\to W_{-1}\to\cdots\to W_{-i}\to U_{-i}\to 0 \\
\tag{$\ast_i$}
0\to M_{-(i+2)}\to U_{-i}\to V_{-i}\to 0\\
\tag{$\circledast_i$}
0\to U_{-i}\to W_{-(i+1)}\to U_{-(i+1)}\to 0 \\
\tag{$\ddag_i$}
0\to M_{-(i+3)}\to U_{-(i+1)}\to V_{-(i+1)}\to 0
\end{gather}
satisfying the following properties:
\begin{enumerate} 
\item[$(\text{a}_i)$] the objects $W_0,\ldots,W_{-(i+1)}$ are in $\catw$;
\item[$(\text{b}_i)$] the sequence $(\dag_i)$ is $\ahom(-,\catw)$-exact;
\item[$(\text{c}_i)$] one has $U_{-i}\perp\catw$ and $U_{-(i+1)}\perp\catw$;
\item[$(\text{d}_i)$] one has $V_{-i}\perp\catw$ and $V_{-(i+1)}\perp\catw$;
\item[$(\text{e}_i)$] one has $V_{-i}$ and $V_{-(i+1)}$ in $\propcorescatw$. 
\end{enumerate}
The sequence $(\dag_i)$ is obtained by splicing the sequences 
$(\dag_0)$, $(\circledast_0)$, \ldots, $(\circledast_{i-1})$.
Continuing to splice inductively,
conditions $(\text{a}_i)-(\text{c}_i)$ show that this process yields
a proper $\catw$-coresolution of $M$, as desired.

We begin with the base case $i=0$. The membership
$G_{-1}\in\propcorescatw$ yields
a proper $\catw$-coresolution of $G_{-1}$
and hence an exact sequence
\begin{equation} \label{exact09}
0\to G_{-1}\to W_0\to V_0\to 0
\end{equation}
that is  $\ahom(-,\catw)$-exact 
and with objects $W_0\in \catw$ and $V_0\in \propcorescatw$. 
Using the conditions $G_{-1}\perp\catw$ and $W_0\perp\catw$,
Lemma~\ref{perp03}\eqref{perp03item1} implies $V_0\perp\catw$.
Consider the pushout diagram whose top row is $(\maltese_{-1})$
and whose middle column is~\eqref{exact09}.
\begin{equation}
\begin{split} \label{diag08}
\xymatrix{
& & 0 \ar[d] & 0 \ar[d] \\
0\ar[r]  & M\ar[r] \ar[d]_{\cong} & G_{-1}\ar[r] \ar[d] \ar@{}[rd]|>>{\lrcorner} & M_{-2}\ar[r]\ar[d]  & 0 \\
0\ar[r] & M\ar[r] & W_{0}\ar[r]\ar[d] & U_0\ar[r]\ar[d] & 0 \\
& & V_0\ar[r]^{\cong}\ar[d] & V_0\ar[d] \\
& & 0 & 0
}
\end{split}
\end{equation}
Applying Lemma~\ref{perp03}\eqref{perp03item1} to the rightmost
column of this diagram,
the conditions $M_{-2}\perp\catw$ and $V_0\perp\catw$ imply
$U_0\perp\catw$.  
For each object $W'\in\catw$, 
use the condition $U_0\perp W'$ with
the long exact sequence in
$\aext(-,W')$ asociated to the middle row of this diagram
to conclude that this row is $\ahom(-,\catw)$-exact.
Set $(\dag_0)$ equal to the middle row of~\eqref{diag08}, and set
$(\ast_0)$ equal to the rightmost
column of~\eqref{diag08}.
Construct the next pushout diagram using $(\maltese_{-2})$ in the top row 
and the rightmost column of~\eqref{diag08} in the left column.
\begin{equation}
\begin{split} \label{diag100}
\xymatrix{
& 0 \ar[d] & 0 \ar[d] \\
0\ar[r]  & M_{-2}\ar[r] \ar[d] \ar@{}[rd]|>>{\lrcorner} & G_{-2}\ar[r] \ar[d] & M_{-3}\ar[r]\ar[d]_{\cong}  & 0 \\
0\ar[r] & U_{0}\ar[r]\ar[d] & Z_0\ar[r]\ar[d] & M_{-3}\ar[r] & 0 \\
& V_{0}\ar[r]^{\cong}\ar[d] & V_{0}\ar[d] \\
& 0 & 0
}
\end{split}
\end{equation}
As in the above discussion,
the condition $V_0\perp\catw$ implies that the middle column of~\eqref{diag100}
is $\ahom(-\catw)$-exact.  
Using Lemma~\ref{perp03}\eqref{perp03item1} with
this column,
the conditions $G_{-2}\perp\catw$ and $V_0\perp\catw$ yield
$Z_0\perp\catw$.
We know that $G_{-2}$ and $V_0$ are in $\propcorescatw$, so 
an application of
Lemma~\ref{horseshoe01}\eqref{horseshoe01item2} to this column implies 
$Z_0\in \propcorescatw$. 
A proper $\catw$-coresolution of $Z_0$
provides  an exact sequence
\begin{equation} \label{exact100}
0\to Z_0\to W_{-1}\to V_{-1}\to 0
\end{equation}
that is  $\ahom(-,\catw)$-exact 
and with objects $W_{-1}\in \catw$ and $V_{-1}\in \propcorescatw$. 
Again using Lemma~\ref{perp03}\eqref{perp03item1},
the conditions $Z_0\perp\catw$ and $W_{-1}\perp\catw$ 
conspire with the $\ahom(-,\catw)$-exactness of~\eqref{exact100} to
imply
$V_{-1}\perp\catw$.
Consider the next pushout diagram whose top row is 
the middle row of~\eqref{diag100}
and whose middle column is~\eqref{exact100}.
\begin{equation}
\begin{split} \label{diag101}
\xymatrix{
& & 0 \ar[d] & 0 \ar[d] \\
0\ar[r]  & U_0\ar[r] \ar[d]_{\cong} & Z_{0}\ar[r] \ar[d] \ar@{}[rd]|>>{\lrcorner} & M_{-3}\ar[r]\ar[d]  & 0 \\
0\ar[r] & U_0\ar[r] & W_{-1}\ar[r]\ar[d] & U_{-1}\ar[r]\ar[d] & 0 \\
& & V_{-1}\ar[r]^{\cong}\ar[d] & V_{-1}\ar[d] \\
& & 0 & 0
}
\end{split}
\end{equation}
Set $(\circledast_0)$ equal to the middle row of this diagram,
and $(\ddag_0)$ equal to the rightmost column.
Thanks to Lemma~\ref{perp03}\eqref{perp03item1},
the conditions $M_{-3}\perp\catw$ and $V_{-1}\perp\catw$
imply $U_{-1}\perp\catw$.  Thus, the conditions
$(\text{a}_0)-(\text{e}_0)$ are satisfied, establishing the base case.

For the induction step, assume that the exact sequences
$(\dagger_i)$, $(\ast_i)$, $(\circledast_i)$, and $(\ddag_i)$ have
been constructed satisfying the conditions $(\text{a}_i)-(\text{e}_i)$.
Note that condition $(\text{c}_i)$ implies that the sequence
$(\circledast_i)$ is $\ahom(-,\catw)$-exact.
Thus, we may splice together the sequences
$(\dagger_i)$ and $(\circledast_i)$ to construct the sequence
$(\dagger_{i+1})$ 
which is exact and $\ahom(-,\catw)$-exact
and such that
$W_0,\ldots,W_{-(i+1)}\in \catw$.
Also, set $(\ast_{i+1})=(\ddag_i)$.

The next pushout diagram has $(\maltese_{i+3})$ in the top row 
and $(\ddag_i)$ in the left column.
\begin{equation}
\begin{split} \label{diag102}
\xymatrix{
& 0 \ar[d] & 0 \ar[d] \\
0\ar[r]  & M_{-(i+3)}\ar[r] \ar[d] \ar@{}[rd]|>>{\lrcorner} 
& G_{-(i+3)}\ar[r] \ar[d] & M_{-(i+4)}\ar[r]\ar[d]_{\cong}  & 0 \\
0\ar[r] & U_{-(i+1)}\ar[r]\ar[d] & Z_{-(i+1)}\ar[r]\ar[d] & M_{-(i+4)}\ar[r] & 0 \\
& V_{-(i+1)}\ar[r]^{\cong}\ar[d] & V_{-(i+1)}\ar[d] \\
& 0 & 0
}
\end{split}
\end{equation}
With the long exact sequence in $\aext(-,-)$,
the condition $V_{-(i+1)}\perp\catw$ implies that the middle column of~\eqref{diag102}
is $\ahom(-,\catw)$-exact.  
Using Lemma~\ref{perp03}\eqref{perp03item1} with
this column,
the conditions $G_{-(i+3)}\perp\catw$ and $V_{-(i+1)}\perp\catw$ yield
$Z_{-(i+1)}\perp\catw$.
As $G_{-(i+3)}$ and $V_{-(i+1)}$ are in $\propcorescatw$, 
apply
Lemma~\ref{horseshoe01}\eqref{horseshoe01item2} to this column to conclude
$Z_{-(i+1)}\in \propcorescatw$. 
A proper $\catw$-coresolution of $Z_{-(i+1)}$
provides  an exact sequence
\begin{equation} \label{exact102}
0\to Z_{-(i+1)}\to W_{-(i+2)}\to V_{-(i+2)}\to 0
\end{equation}
that is  $\ahom(-,\catw)$-exact 
and with objects $W_{-(i+2)}\in \catw$ and $V_{-(i+2)}\in \propcorescatw$. 
Again using Lemma~\ref{perp03}\eqref{perp03item1} and the
$\ahom(-,\catw)$-exactness of~\eqref{exact102},
the conditions $Z_{-(i+1)}\perp\catw$ and $W_{-(i+2)}\perp\catw$ imply
$V_{-(i+2)}\perp\catw$.
Consider the next pushout diagram whose top row is 
the middle row of~\eqref{diag102}
and whose middle column is~\eqref{exact102}. 
\begin{equation}
\begin{split} \label{diag103}
\xymatrix{
& & 0 \ar[d] & 0 \ar[d] \\
0\ar[r]  & U_{-(i+1)}\ar[r] \ar[d]_{\cong} & Z_{-(i+1)}\ar[r] \ar[d] \ar@{}[rd]|>>{\lrcorner} 
& M_{-(i+4)}\ar[r]\ar[d]  & 0 \\
0\ar[r] & U_{-(i+1)}\ar[r] & W_{-(i+2)}\ar[r]\ar[d] & U_{-(i+2)}\ar[r]\ar[d] & 0 \\
& & V_{-(i+2)}\ar[r]^{\cong}\ar[d] & V_{-(i+2)}\ar[d] \\
& & 0 & 0
}
\end{split}
\end{equation}
Set $(\circledast_{i+1})$ equal to the middle row of this diagram,
and set $(\ddag_{i+1})$ equal to the rightmost column.
With Lemma~\ref{perp03}\eqref{perp03item1},
the conditions $M_{-(i+4)}\perp\catw$ and $V_{-(i+2)}\perp\catw$
imply $U_{-(i+2)}\perp\catw$.  Thus, the conditions
$(\text{a}_{i+1})-(\text{e}_{i+1})$ are satisfied, establishing the induction step.
\end{proof}

What follows is the second step in the proof of Theorem~\ref{thma} from the introduction.
See Example~\ref{Fc} for the necessity of the cogeneration hypothesis.

\begin{thm} \label{gxgw01}
If $\catx$ is closed under extensions
and $\catw$ is both an injective cogenerator and
a projective generator for $\catx$, then
$\opg^n(\catx) \subseteq\opg(\catw)$ for each $n\geq 1$.
\end{thm}

\begin{proof}
Let $N$ be an object in $\opg^n(\catx)$.  
By Theorem~\ref{ggisg01} we know that $N$ 
admits a proper $\catw$-resolution $W'\xra{\simeq}N$
and a proper $\catw$-coresolution $N\xra{\simeq}W''$.
We will show that $^+W'$ is $\ahom(-,\catw)$-exact and
that $^+W''$ is $\ahom(\catw,-)$-exact.
Since we already know that
$^+W'$ is $\ahom(\catw,-)$-exact and
that $^+W''$ is $\ahom(-,\catw)$-exact,
this will show that the concatenated complex
$$\cdots\to W'_1\to W'_0\to W''_0\to W'_{-1}\to\cdots$$
is a complete $\catw$-resolution of $N$, completing the proof.

We will show that $^+W'$ is $\ahom(-,\catw)$-exact; the proof of the other
fact is dual.  For each $i\geq 0$, there is an exact sequence
\begin{equation} \label{exact12} \tag{$\ast_i$}
0\to N_{i+1}\to W_i'\to N_i\to 0
\end{equation}
We have $N_0=N$ and so $N_0\perp\catw$ is true by Proposition~\ref{gperp01};  
and $W'_i\perp\catw$ by assumption.  Using Lemma~\ref{perp03}\eqref{perp03item1}, 
an induction
argument implies $N_i\perp\catw$ for each $i$.  It follows that
($\ast_i$) is $\ahom(-,\catw)$-exact, and it follows that $^+W'$ is $\ahom(-,\catw)$-exact,
as desired.
\end{proof}

Theorem~\ref{thma} from the introduction follows from the next result; see Example~\ref{GP}.

\begin{cor} \label{gxgw02}
If $\catw\perp\catw$, then 
$\opg^n(\catw)=\opg(\catw)$ for each $n\geq 1$.
\end{cor}

\begin{proof}
Note  that Corollaries~\ref{exten02} and~\ref{wcogengw01} imply that
$\opg(\catw)$ is  closed under extensions
and that $\catw$ is both an injective cogenerator and
a projective generator for $\opg(\catw)$.

We argue by induction on $n$, the case 
$n=1$ being trivial.
For $n> 1$, if $\opg^{n-1}(\catw)=\opg(\catw)$, then 
setting $\catx=\opg(\catw)$ in Theorem~\ref{gxgw01} yields the
final containment in the next sequence
$$\opg(\catw)\subseteq\opg^{n}(\catw)=\opg(\opg^{n-1}(\catw))
=\opg(\opg(\catw))\subseteq\opg(\catw)
$$
and hence the desired conclusion.
\end{proof}

With Corollary~\ref{exten02},
the final two results of this section contain Theorem~\ref{thmb} from the introduction;
compare to~\cite[(4.3.5)]{christensen:gd}

\begin{prop} \label{exact01}
If $\catw\perp\catw$, then 
$\opg(\catw)$ is closed under direct summands.
\end{prop}

\begin{proof}
Let $A'$ and $A''$ be objects in $\cata$ such that $A'\oplus A''$ is in $\opg(\catw)$.
We construct
proper $\catw$-resolutions $W'\xra{\simeq}A'$ and $W''\xra{\simeq}A''$
such that $(W')^+$ and $(W'')^+$ are $\ahom(-,\catw)$-exact.
Dually, one constructs proper $\catw$-coresolutions $A'\xra{\simeq}V'$ and $A''\xra{\simeq}V''$
such that $^+V'$ and $^+V''$ are 
$\ahom(\catw,-)$-exact, and this shows that $A'$ and $A''$ are in $\opg(\catw)$.  

Observe first that $A'$ and $A''$ both admit (augmented) proper $\opg(\catw)$-resolutions
\begin{align*}
X'&=
\cdots\xra{\left(\begin{smallmatrix} \id & 0 \\ 0 & 0 \end{smallmatrix}\right)}
A'\oplus A''\xra{\left(\begin{smallmatrix} 0 & 0 \\ 0 & \id \end{smallmatrix}\right)}
A'\oplus A''\xra{\left(\begin{smallmatrix} \id & 0 \end{smallmatrix}\right)} A'\to 0
\\
X''&=
\cdots\xra{\left(\begin{smallmatrix} 0 & 0 \\ 0 & \id \end{smallmatrix}\right)}
A'\oplus A''\xra{\left(\begin{smallmatrix} \id & 0 \\ 0 & 0 \end{smallmatrix}\right)}
A'\oplus A''\xra{\left(\begin{smallmatrix} 0  &  \id \end{smallmatrix}\right)}A''\to 0
\end{align*}
where properness follows from the contractibility of $X'$ and $X''$.
From Proposition~\ref{gperp01} we know 
$\catw\perp (A'\oplus A'')$, so the additivity of $\aext$ implies
$\catw\perp A'$ and
$\catw\perp A''$.
Lemma ~\ref{coresxw01}\eqref{coresxw01item2}
and Corollary~\ref{wcogengw01} imply
that $A'$ and $A''$ admit proper $\catw$-resolutions
$W'\xra{\simeq}A'$ and $W''\xra{\simeq}A''$, so
$W'\oplus W''\xra{\simeq}A'\oplus A''$
is a proper $\catw$-resolution.
We show that
$(W')^+$ and $(W'')^+$ are $\ahom(-,\catw)$-exact.
As 
$A'\oplus A''$ is in $\opg(\catw)$, it admits a proper
$\catw$-resolution $W\xra{\simeq}A'\oplus A''$
such that $W^+$ is $\ahom(-,\catw)$-exact.
Hence, the resolutions $W$ and $W'\oplus W''$ are homotopy equivalent.
Because $W^+$ is $\ahom(-,\catw)$-exact, we know that
$(W'\oplus W'')^+$ is also $\ahom(-,\catw)$-exact,
and so $(W')^+$ and $(W'')^+$ are $\ahom(-,\catw)$-exact.
\end{proof}

\begin{thm} \label{admepi01}
Assume $\catw\perp\catw$.
\begin{enumerate}[\quad\rm(a)]
\item \label{admepi01item1}
If $\catw$ is closed under kernels of epimorphisms, then so is
$\opg(\catw)$.
\item \label{admepi01item2}
If $\catw$ is closed under cokernels of monomorphisms, then so is
$\opg(\catw)$.
\end{enumerate}
\end{thm}

\begin{proof}
We prove part~\eqref{admepi01item1}; the proof of part~\eqref{admepi01item2} is dual.
Consider an exact sequence in $\cata$ with objects $N$ and $N''$ in $\opg(\catw)$.
\begin{equation} \label{exact16}
0\to N'\to N\xra{\tau} N''\to 0
\end{equation}
Let $W$ and $W''$ be complete $\catw$-resolutions of $N$ and $N''$, respectively.

We first construct a commutative diagram of the following form
\begin{equation} \label{diag09}
\begin{split}
\xymatrix{
\cdots \ar[r]^{\partial^{W}_{2}} & W_1 \ar[r]^{\partial^{W}_{1}}\ar@{-->}[d]_{\wti{\tau}_{1}} 
& W_0\ar[r]^{\pi}\ar@{-->}[d]_{\wti{\tau}_{0}} & N \ar[r]^{\epsilon}\ar[d]_{\tau} 
& W_{-1} \ar[r]^{\partial^{W}_{-1}}\ar@{-->}[d]_{\wti{\tau}_{-1}} 
& W_{-2}\ar[r]^{\partial^{W}_{-2}}\ar@{-->}[d]_{\wti{\tau}_{-2}} 
& \cdots \\
\cdots \ar[r]^{\partial^{W''}_{2}} & W''_1 \ar[r]^{\partial^{W''}_{1}} & W''_0\ar[r]^{\pi''} & N'' \ar[r]^{\epsilon''}
& W''_{-1} \ar[r]^{\partial^{W''}_{-1}} & W''_{-2}\ar[r]^{\partial^{W''}_{-2}} & \cdots 
}
\end{split}
\end{equation}
where $\epsilon\pi=\partial^W_0$ and $\epsilon''\pi''=\partial^{W''}_0$.
Since $(W_{\geq 0})^+$ is a chain complex and $(W''_{\geq 0})^+$ is $\ahom(\catw,-)$ exact,
one can successively lift $\tau$ to the left as in the diagram; argue as in~\cite[(1.8)]{holm:ghd}.  
Dually, since $^+W''_{< 0}$ is a chain complex and $^+W_{< 0}$ is $\ahom(-,\catw)$ exact,
one can successively lift $\tau$ to the right as in the diagram.  

Thus, we have constructed a morphism of chain complexes $\wti{\tau}\colon W\to W''$
such that the induced map $\coker(\partial^{W}_1)\to\coker(\partial^{W''}_1)$ is equivalent
to $\tau$.

Next, we show that
there exists a complex $\wti{W}$ with
and a morphism $\tau'\colon W\oplus \wti{W}\to W''$
satisfying the following properties:
\begin{enumerate}[\quad(a)]
\item \label{admepi011}
$\wti{W}$ is contractible and $\wti{W}_n$ is  in $\catw$  for each $n\in\mathbb{Z}$.
\item \label{admepi012}
$\tau_n'$ is an epimorphism for each $n\in\mathbb{Z}$.
\item \label{admepi013}
The natural monomorphism $W\xra{\epsilon} W\oplus \wti{W}$ 
satisfies $\wti{\tau}=\tau'\epsilon$.
\end{enumerate}
The complex 
$\wti{W}=\shift^{-1}\cone(\id_{W''})$ is contractible; see~\ref{notation07a}.
Let $f\colon \wti{W}\to W''$ denote the composition of the natural morphisms 
$\wti{W}=\shift^{-1}\cone(\id_{W''})\to W''$.
Note that each $f_n$ is a split epimorphism.
It follows that the homomorphisms
$\tau'_n=(\wti{\tau}_n \,\, f_n)\colon W_n\oplus \wti{W}_n\to W''_n$
describe a morphism of complexes satisfying the desired properties.

Because of property~\eqref{admepi011} 
the complex $\wti{W}$ is a complete $\catw$-resolution; see Remark~\ref{gtriv01}.
Set $\wti{N}=\coker(\partial_1^{\wti{W}})$, which is an object in $\opg(\catw)$
with complete resolution $\wti{W}$.
Furthermore, one has $\coker(\partial_1^{W\oplus \wti{W}})\cong N\oplus\wti{N}$,
and the morphism $\tau'$ induces a homomorphism
$N\oplus\wti{N}\xra{f=(\tau\,\,\pi)}N''$.
Because $\tau$ is surjective, the map $f$ is also surjective.
We will show that $\Ker(f)$ is in $\opg(\catw)$, and then we will show that
$N'=\ker(\tau)$ is in $\opg(\catw)$.

The morphism $\tau'$ is degreewise surjective.  
As $\catw$ is closed under kernels of epimorphisms, it follows
that the complex $W'=\ker(\wti{\tau})$ consists of objects in $\catw$.
The next exact sequence shows that $W'$ is exact because
$W$, $\wti{W}$, and $W''$ are so
\begin{equation*} 
0\to W'\to W\oplus \wti{W}\xra{\tau'} W''\to 0.
\end{equation*}
This sequence  induces a second exact sequence
$$0\to W'_{\geq 0}\to W_{\geq 0}\oplus \wti{W}_{\geq 0}\to W''_{\geq 0}\to 0$$
whose associated
long exact sequence is of the form
$$0\to \coker(\partial^{W'}_1)\to N\xra{f}N''\to 0.$$
Thus, we have $\Ker(f)\cong \coker(\partial^{W}_1)$.
To show that $\Ker(f)$ is in $\opg(\catw)$,
it suffices to show that  $W'$ is 
$\ahom(\catw,-)$-exact and $\ahom(-,\catw)$-exact. 
For each object $U$ in $\catw$, the next sequence
of complexes is exact as $\catw\perp\catw$
$$0\to \ahom(U,W')\to \ahom(U,W)\to \ahom(U,W'')\to 0.$$
Since $W$ and $W''$ are $\ahom(\catw,-)$-exact, the associated 
long exact sequence shows that 
$W'$ is also $\ahom(\catw,-)$-exact.
Dually, one shows that $W'$ is $\ahom(-,\catw)$-exact, thus 
showing that $\Ker(f)$ is in $\opg(\catw)$.

To see that $N'$ is in $\opg(\catw)$, consider the following pullback diagram 
whose rightmost column is~\eqref{exact16} and whose middle row is the
natural split exact sequence.
\begin{equation} \label{diag10}
\begin{split}
\xymatrix{
& & 0 \ar[d] & 0 \ar[d] \\
0\ar[r] & \wti{N}\ar[r]^-{\gamma}\ar[d]_{\cong}^{\alpha} & \Ker(f)\ar[r]\ar[d]^{\delta} \ar@{}[rd]|<<{\ulcorner} & N'\ar[r]\ar[d] & 0 \\
0\ar[r] & \wti{N}\ar[r]^-{\beta} & N\oplus\wti{N} \ar[r]\ar[d]^f & N \ar[r]\ar[d]^{\tau} & 0 \\
&& N'' \ar[r]^{\cong} \ar[d] & N'' \ar[d] \\
& & 0 & 0
}
\end{split}
\end{equation}
Let $\sigma\colon N\oplus\wti{N}\to \wti{N}$ denote the natural surjection.
It follows that $\sigma\beta=\id_{\wti{N}}$.
Since $\alpha$ is an isomorphism, the equality $\delta\gamma=\beta\alpha$
implies 
$$(\alpha^{-1} \sigma \delta) \gamma=\alpha^{-1} \sigma\beta\alpha
=\alpha^{-1}\alpha=\id_{\wti{N}}$$
and so the top row of~\eqref{diag10} is split exact.
Hence, the object $N'$ is a direct summand of $\Ker(f)$.
We have shown that  $\Ker(f)$ is in $\opg(\catw)$.
The category $\opg(\catw)$ is closed under direct summands by Proposition~\ref{exact01},
and so $N'$ is in $\opg(\catw)$ as desired.
\end{proof}

\section{Consequences for categories of interest} \label{sec5}

Let $R$ be a
commutative ring and $C$ a semidualizing $R$-module.
We now apply the results of Section~\ref{sec3} 
to the examples in Section~\ref{sec4}.  We begin
with some computations.

\begin{ex} \label{gp01}
The relevant definitions yield equalities $\opg(\catp(R))=\catgp(R)$
and $\opg(\cati(R))=\catgi(R)$.  If $R$ is noetherian,
then $\opg(\catp^f(R))=\catg(R)$.
\end{ex}

The next result generalizes the previous example.

\begin{prop} \label{gpc01}
Let $R$ be a commutative ring.
If   $C$ is $R$-semidualizing, then
$\opg(\catpc(R))=\catgpc(R)\cap\catbc(R)$ and
$\opg(\catic(R))=\catgic(R)\cap\catac(R)$.
If further $R$ is noetherian, then
$\opg(\catpc^f(R))=\catg_C(R)\cap\catbc(R)$.
\end{prop}

\begin{proof}
We will prove the first equality;  the others are proved similarly.  For one containment, let
$M$ be an object in $\opg(\catpc(R))$.  
To show that $M$ is an object in $\catgpc(R)$, we use~\cite[(3.2)]{white:gpdrsm}:
it suffices to show that $M$ admits a proper $\catpc(R)$-coresolution and
$M\perp\catpc(R)$.  The first of these is in Remark~\ref{gtriv01}
which says that $M$ is in $\propcorescatpc$; the second
one is in Proposition~\ref{gperp01} which implies $\opg(\catpc)\perp\catpc$.
To show that $M$ is an object in $\catbc(R)$, we need to verify
$\ext_R^{\geq 1}(C,M)=0$ and
$\tor_{\geq 1}^R(C,\Hom_R(C,M))=0$
and $M\cong C\otimes_R\Hom_R(C,M)$.
The first of these is in Proposition~\ref{gperp01} which implies $\catpc\perp\opg(\catpc)$,
and the others are in~\cite[(2.2)]{takahashi:hasm}.

For the reverse containment, fix an object $N$ in $\catgpc(R)\cap\catbc(R)$.
Since $N$ is in $\catgpc(R)$, it admits
a complete $\catp\catpc$-resolution  $Y$, so the complex $Y_{<0}$ is a proper
$\catpc(R)$-coresolution of $N$.
Also, $N$ admits a proper $\catpc(R)$-resolution $Z$ by~\cite[(2.4)]{takahashi:hasm}
as $N$ is in $\catbc(R)$.
Once it is shown that $Y_{<0}^+$ is $\Hom_R(\catpc(R),-)$-exact 
and $Z^+$ is $\Hom_R(-,\catpc(R))$-exact, a complete $\catp\catpc$-resolution of $N$
will be obtained by splicing $Z$ and $Y_{<0}$.

To see that $Y_{<0}^+$ is $\Hom_R(\catpc(R),-)$-exact, set $N^{(0)}=N$ and 
$N^{(i)}=\coker(\partial^Y_{i-1})$ for each $i\leq -1$.  
From~\cite[(5.2)]{holm:fear}, we know that  $Y_i$ is in $\catbc(R)$ for each $i\leq -1$.
Since $N$ is also in $\catbc(R)$,
an induction argument using the exact sequence
\begin{equation} \label{exact101} \tag{$\ast_i$}
0\to N^{(i-1)}\to Y_{i}\to N^{(i)}\to 0
\end{equation}
implies that $N^{(i)}$ is in $\catbc(R)$ for each $i\leq -1$.
For each projective $R$-module $P$, this yields the vanishing in the next sequence
$$\ext^1_R(P\otimes_R C,N^{(i)})\cong\Hom_R(P,\ext^1_R(C,N^{(i)}))=0$$
while the isomorphism is from Hom-Tensor adjunction.  It follows that
each sequence~\eqref{exact101} is $\Hom_R(\catpc(R),-)$-exact, and
thus so is  $Y_{<0}^+$.

To see that $Z^+$ is $\Hom_R(-,\catpc(R))$-exact, 
it suffices to let $P$ be projective and 
to justify the following sequence for $i\geq 1$.
$$\HH_{-i}(\Hom_R(Z^+,P\otimes_R C))=\HH_{-i}(\Hom_R(Z,P\otimes_R C))
\cong\ext^i_R(N,P\otimes_R C)
=0$$
The isomorphism is from~\cite[(4.2.a)]{takahashi:hasm} because
$N$ and $P\otimes_R C$ are in $\catbc(R)$.  The vanishing follows because
$N$ is in $\catgpc(R)$
and $\catgpc(R)\perp\catpc(R)$; see~\cite[(3.2)]{white:gpdrsm}.
\end{proof}

We now outline the consequences of Corollaries~\ref{wcogengw01}
and~\ref{gxgw02} for the examples of Section~\ref{sec4}.  The first example
below
contains Theorem~\ref{thma} from the introduction.

\begin{ex} \label{GP}
The category $\catpc(R)$ is an injective cogenerator and
a projective generator for $\catgpc(R)\cap\catbc(R)$, and
$\opg^n(\catpc(R))=\catgpc(R)\cap\catbc(R)$ for each $n\geq 1$.
Hence, $\catp(R)$ is an injective cogenerator and
a projective generator for $\catgp(R)$, and  $\opg^n(\catp(R))=\catgp(R)$.  
If $\cata$ has enough projectives, then $\catp(\cata)$
is an injective cogenerator and
a projective generator for $\opg(\catp(\cata))$, and $\opg^n(\catp(\cata))=\opg(\catp(\cata))$.
\end{ex}

\begin{ex} \label{GI}
The category $\catic(R)$ is an injective cogenerator and
a projective generator for $\catgic(R)\cap\catac(R)$, and 
$\opg^n(\catic(R))=\catgic(R)\cap\catac(R)$ for each $n\geq 1$.
Hence, $\cati(R)$ is an injective cogenerator and
a projective generator for $\catgi(R)$, and $\opg^n(\cati(R))=\catgi(R)$.
If $\cata$ has enough injectives, then $\cati(\cata)$
is an injective cogenerator and
a projective generator for $\opg(\cati(\cata))$,
and $\opg^n(\cati(\cata))=\opg(\cati(\cata))$.
\end{ex}

\begin{ex} \label{GC}
Assume that $R$ is noetherian. Then
$\catpc^f(R)$ is an injective cogenerator and
a projective generator for $\catg_C(R)\cap\catbc(R)$,
and $\opg^n(\catpc^f(R))=\catg_C(R)\cap\catbc(R)$.
In particular, $\catp^f(R)$ is an injective cogenerator and
a projective generator for $\catg(R)$, and $\opg^n(\catp^f(R))=\catg(R)$.
\end{ex}

\begin{ex} \label{GFc}
The category $\catfc(R)$ is an injective cogenerator and
a projective generator for $\catg(\catfc(R))$, and 
$\opg^n(\catfc(R))=\catg(\catfc(R))$ for each $n\geq 1$.
\end{ex}

With Proposition~\ref{gperp01} and Corollary~\ref{exten02} in mind,
we now show that $\catw\perp\catw$ need not imply
$\opg(\catw)\perp\opg(\catw)$.

\begin{ex} \label{gexten01}
Let $(R,\m,k)$ be a local, nonregular, Gorenstein, artinian ring. 
With $\catw=\catp^f(R)$, we have $\opg(\catw)=\catg(R)=\catm^f(R)$
where the last equality is because $R$ is artinian and Gorenstein;
see~\cite[(1.4.8),(1.4.9)]{christensen:gd}.
We know $\ext_R^{\geq 1}(k,k)\neq 0$ since $R$ is nonregular, and so
$\opg(\catw)\not\perp\opg(\catw)$.
\end{ex}

We conclude with some questions and final observations.

\begin{question} \label{stab01}
Must there be an equality $\gnx{n}=\gx$ for each $n\geq 1$?
Is $\gx$ always exact?
Is $\gx$ always closed under kernels of epimorphisms
or cokernels of monomorphisms?
Must $\gw$ be contained in $\gx$?
Can  $\catg(\catf(R))$
or $\catg(\catfc(R))$
or $\catg(\catg'(R))$ be identified
as in Proposition~\ref{gpc01}? 
\end{question}

The final examples are presented with an eye toward the
last question in~\ref{stab01}.

\begin{ex} \label{Fc}
If $(R,\fm)$ is a noetherian local ring
and $\dim(R)\geq 1$, then $\catg(\catf(R))\not\subseteq \catg(\catp(R))$. 
Indeed, the ring of formal power series $R[\![X]\!]$ is a flat $R$-module, so
it is in $\catg(\catf(R))$.
Suppose by way of contradiction 
that $R[\![X]\!]$ is in $\catg(\catp(R))$. 
First note that~\cite[Prop.~6]{jensen:vl}
yields $\pdim_R(R[\![X]\!])<\infty$, and so~\cite[(2.7)]{holm:ghd} implies
$\pd_R(R[\![X]\!])=\gpd_R(R[\![X]\!])=0$.  It follows that $R[\![X]\!]$ is projective,  
contradicting~\cite[(2.1)]{buchweitz:psrp}.  

From this it follows that 
the conclusion of Theorem~\ref{gxgw01} need not hold if $\catw$ is not a cogenerator
for $\catx$.
To see this, assume that $R$ is $\fm$-adically complete.  
Standard results combine to show that $\catp(R)$ is a projective generator for
$\catf(R)$ and that $\catf(R)$ is closed under extensions.  Furthermore,
one has $\catf(R)\perp\catp(R)$ by~\cite[(5.3.28)]{enochs:rha}. 
\end{ex}

With Theorem~\ref{gxgw01}, the previous example provides the next result.

\begin{cor} \label{notcogen}
If $R$ is a complete local notherian ring and $\dim(R)\geq 1$, then
$\catp(R)$ is not a cogenerator for
$\catf(R)$. \qed
\end{cor}

\begin{ex} \label{Fc2}
Let $(R,\fm)$ be a noetherian local ring.
If $R$ is not $\fm$-adically complete, then $\catg(\catfc(R))\not\subseteq \catg(\catp(R))$.
The $\fm$-adic completion $\comp{R}$ is flat and cotorsion; see, e.g., \cite[(5.3.28)]{enochs:rha}.
Arguing as in Example~\ref{Fc}, it then suffices to note that $\comp{R}$ is not projective
by~\cite[Thm.~A]{frankild:dcev}.
\end{ex}

\begin{ex} \label{lowerci01}
Assume that $R$ is 
local and noetherian.
Using Example~\ref{GP} and~\cite[(2.8)]{gerko:ohd}, it is straightforward to show that
$\catp^f(R)$ is an injective cogenerator 
and a projective generator for $\catg'(R)$ and that $\catg'(R)$
is closed under extensions.
Theorem ~\ref{gxgw01} now implies 
$\opg^n(\catg'(R)) \subseteq\catg(R)$ for each $n\geq 1$.
\end{ex}

\section*{Acknowledgments}

We are indebted to the referee for 
his/her careful reading of this work.


\providecommand{\bysame}{\leavevmode\hbox to3em{\hrulefill}\thinspace}
\providecommand{\MR}{\relax\ifhmode\unskip\space\fi MR }
\providecommand{\MRhref}[2]{%
  \href{http://www.ams.org/mathscinet-getitem?mr=#1}{#2}
}
\providecommand{\href}[2]{#2}

\end{document}